# Alternatives To Pearson's and Spearman's Correlation Coefficients

Florentin Smarandache
Chair of Math & Sciences Department
University of New Mexico
Gallup, NM 87301, USA

**Abstract**. This article presents several alternatives to Pearson's correlation coefficient and many examples. In the samples where the rank in a discrete variable counts more than the variable values, the mixture of Pearson's and Spearman's gives a better result.

## Introduction

Let's consider a bivariate sample, which consists of n ≥ 2 pairs (x,y). We denote these pairs by:
$(x_1, y_1), (x_2, y_2), \ldots , (x_n, y_n)$,

where $x_i$ = the value of x for the *i*-th observation,
and $y_i$ = the value of y for the *i*-th observation,
for any $1 \leq i \leq n$.

We can construct a scatter plot in order to detect any relationship between variables x and y, drawing a horizontal x-axis and a vertical y-axis, and plotting points of coordinates $(x_i, y_i)$ for all $i \in \{1, 2, \ldots, n\}$.

We use the standard statistics notations, mostly used in regression analysis:

$$\sum x = \sum_{i=1}^{n} x_i, \quad \sum y = \sum_{i=1}^{n} y_i, \quad \sum xy = \sum_{i=1}^{n} (x_i y_i),$$

$$\sum x^2 = \sum_{i=1}^{n} x_i^2, \quad \sum y^2 = \sum_{i=1}^{n} y_i^2, \qquad (1)$$

$$\bar{X} = \frac{\sum_{i=1}^{n} x_i}{n} = \text{the mean of sample variable x},$$

$$\bar{Y} = \frac{\sum_{i=1}^{n} y_i}{n} = \text{the mean of sample variable y}.$$

Let's introduce a notation for the median:

$X_M$ = the median of sample variable x, (2)

$Y_M$ = the median of sample variable y.

**<u>Correlation Coefficients.</u>**

Correlation coefficient of variables x and y shows how strongly the values of these variables are related to one another. It is denoted by r and $r \in [-1, 1]$.

If the correlation coefficient is positive, then both variables are simultaneously increasing (or simultaneously decreasing).

If the correlation coefficient is negative, then when one variable increases while the other decreases, and reciprocally.

Therefore, the correlation coefficient measures the degree of line association between two variables.

We have <u>strong relationship</u> if $r \in [0.8, 1]$ or $r \in [-1, -0.8]$;
<u>moderate relationship</u> if $r \in (0.5, 0.8)$ or $r \in (-0.8, -0.5)$; (3)
And <u>weak relationship</u> if $r \in [-0.5, 0.5]$.

Correlation coefficient does not depend on the measurement unit, neither on the order of variables: (x, y) or (y, x).

If r = 1 or -1, then there is a perfectly linear relationship between x and y. If r = 0, or close to zero, then there is not a strong <u>linear</u> relationship, but there might be a strong <u>non-linear</u> relationship that can be checked on the scatter plot.

The <u>coefficient of determination</u>, denoted by $r^2$, represents the proportion of variation in *y* due to a linear relationship between x and y in the sample:

$$r^2 = \frac{SSTo - SS\,Resid}{SSTo} = 1 - \frac{SSResid}{SSTo} \qquad (4)$$

where SSTo = total sum of squares = $\sum (y - \bar{y})^2 = \sum_{i=1}^{n}(y_i - \bar{y})^2$ (5)

and SSResid = residual sum of squares = $\sum (y - \hat{y})^2 = \sum_{i=1}^{n}(y_i - \hat{y}_i)$ (6)

with $\hat{y}_i$ = the i-th predicted value = $a + bx_i$ for $i \in \{1, 2, \dots, n\}$

resulting from substituting each sample x value into the equation for the least-squares line

ŷ = a + bx

$$\text{where } b = \frac{\sum xy - [(\sum x)(\sum y)/n]}{\sum x^2 - [(\sum x)^2/n]} \qquad (7)$$

$$\text{and } a = \bar{Y} - b\bar{X}. \qquad (8)$$

Obviously: coefficient of determination = (correlation coefficient)$^2$.

Two sample correlation coefficients are well-known:

1) <u>Pearson's sample correlation coefficient</u>, let's denote it by $r_p$

$$r_p = \frac{\sum xy - [(\sum x)(\sum y)/n]}{\sqrt{\sum x^2 - [(\sum x)^2/n]} \cdot \sqrt{\sum y^2 - [(\sum y)^2/n]}} \qquad (9)$$

which is the most popular;

and 2) <u>Spearman's rank correlation coefficient</u>, let's denote it by $r_S$, which is obtained from the previous one by replacing, for each $i \in \{1, 2, ..., n\}$, $x_i$ by its rank in the variable x, and similarly for $y_i$.

\*

We propose more <u>alternative sample correlation coefficients</u> in the following ways, replacing in Pearson's formula (9):

3.1. Each $x_i$ by its deviation from the x mean: $x_i - \bar{x}$,
and each $y_i$ by its deviation from the y mean: $y_i - \bar{y}$.

3.2. Each $x_i$ by its deviation from the x minimum: $x_i - x_{min}$, and each $y_i$ by its deviation from the y minimum: $y_i - y_{min}$.

3.3. Each $x_i$ by its deviation from the x maximum: $x_{max} - x_i$, and each $y_i$ by its deviation from the y maximum: $y_{max} - y_i$

3.4. Each $x_i$ by its deviation from a given $x_k$ (for $k \in \{1, 2, ..., n\}$):

$$x_i - x_k$$
and each $y_i$ by its deviation from the corresponding given $y_k$:
$$y_i - y_k$$

Not surprisingly, all these four alternative sample correlation coefficients are equal to Pearson's since they are simply related to translations of Cartesian axes, whose origin (0,0) is moved to $(\bar{x}, \bar{y})$, $(x_{min}, y_{min})$, $(x_{max}, y_{max})$, or $(x_k, y_k)$ respectively.

**Example**: Let the variables x, y be given below:

| x | 6 | 7 | 12 | 14 | 23 | 41 | 53 | 60 | 69 | 72 |
|---|---|---|----|----|----|----|----|----|----|----|
| y | 2.5 | 1.1 | 6.3 | 2.1 | 2.9 | 15.3 | 20.7 | 18.4 | 22 | 33 |

Table 1

and their scatter plot:

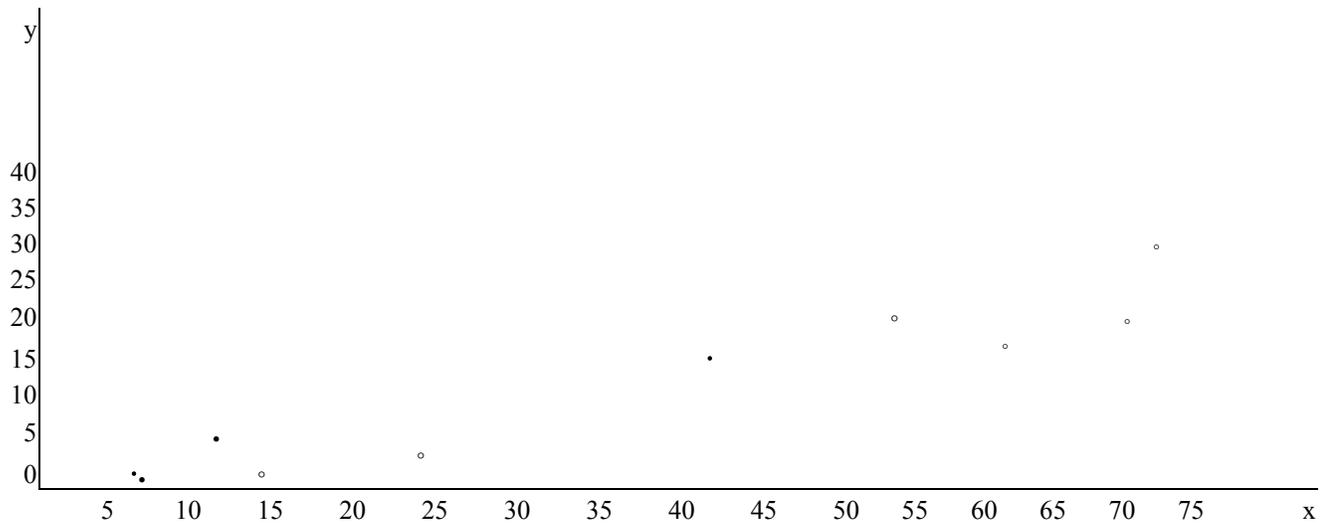

Graph 1

1) Calculating <u>Pearson's</u> correlation coefficient:

$$\sum x = 357; \qquad \bar{x} = 35.7;$$

$$\sum y = 124.3; \qquad \bar{y} = 12.43;$$

$$\sum x^2 = 18,989;$$

$$\sum y^2 = 2,634.11;$$

$$\sum xy = 6,916.8;$$

$r_p = 0.95075$.

2) Calculating <u>Spearman's</u> rank correlation coefficient:

| x | 1 | 2 | 3 | 4 | 5 | 6 | 7 | 8 | 9 | 10 |
|---|---|---|---|---|---|---|---|---|---|----|
| y | 3 | 1 | 5 | 2 | 4 | 6 | 8 | 7 | 9 | 10 |

Table 2

$$\sum x = \frac{(1+10) \cdot 10}{2} = 11.5 = 5.5;$$

$$\sum y = 55;$$

$$\sum x^2 = 385;$$

$$\sum y^2 = 385;$$

$$\sum xy = 377;$$

$r_s = 0.90303$.

3.1) Replacing $x_i$ by $x_i - \bar{x}$ and $y_i$ by $y_i - \bar{y}$ for all $i$ (deviations from the mean):

| x | -29.7 | -28.7 | -23.7 | -21.7 | -12.7 | 5.3 | 17.3 | 24.3 | 33.3 | 36.3 |
|---|-------|-------|-------|-------|-------|-----|------|------|------|------|
| y | -9.93 | -11.33 | -6.13 | -10.33 | -9.53 | 2.87 | 8.27 | 5.97 | 9.57 | 20.57 |

Table 3

Similarly: $\sum x = 0$,

because $\sum x = \sum_{i=1}^{10}(x_i - \bar{x}) = x_1 - \bar{x} + x_2 - \bar{x} + \ldots + x_{10} - \bar{x} = (x_1 + x_2 + \ldots + x_{10}) - 10\bar{x}$

$= (x_1 + x_2 + \ldots + x_{10}) - 10 \cdot \frac{x_1 + x_2 + \ldots + x_n}{10} = 0;$

$\sum y = 0;$
$\sum x^2 = 6,244.10;$
$\sum y^2 = 1,089.06;$

$$\sum xy = 2{,}479.29;$$

$$r_{mean} = 0.95075.$$

3.2) Replacing $x_i$, $y_i$ by their deviations from the smaller x: = x-$x_{small}$ and y: = y-$y_{small}$ we have a translation of axes again.

| x | 0 | 1 | 6 | 8 | 17 | 35 | 47 | 54 | 63 | 66 |
|---|---|---|---|---|----|----|----|----|----|----|
| y | 1.4 | 0 | 5.2 | 1 | 1.8 | 14.2 | 19.6 | 17.3 | 20.9 | 31.9 |

Table 4

$$\sum x = 297;$$
$$\sum y = 113.3;$$
$$\sum x^2 = 15{,}065;$$
$$\sum y^2 = 2{,}372.75;$$
$$\sum xy = 5{,}844.30;$$
$$r_{(small)} = 0.95075.$$

3.3) Replacing $x_i$, $y_i$ by their deviations from the maximum:

| x | 66 | 65 | 60 | 58 | 49 | 31 | 19 | 12 | 3 | 0 |
|---|----|----|----|----|----|----|----|----|---|---|
| y | 30.5 | 31.9 | 26.7 | 30.9 | 30.1 | 17.7 | 12.3 | 14.6 | 11 | 0 |

Table 5

$$\sum x = 363;$$
$$\sum y = 205.7;$$
$$\sum x^2 = 19{,}421;$$
$$\sum y^2 = 5{,}320.31;$$
$$\sum xy = 9{,}946.20;$$
$$r_{(max)} = 0.95075.$$

3.4) Replacing $x_i$ by $x_i - x_4$ and $y_i$ by $y_i - y_4$ (in this case k = 4), $(x_4, y_4) = (14, 2.1)$:

| x | -8 | -7 | -2 | 0 | 9 | 27 | 39 | 46 | 55 | 58 |
|---|----|----|----|---|---|----|----|----|----|----|
| y | 0.4 | -1 | 4.2 | 0 | 0.8 | 13.2 | 18.6 | 16.3 | 19.9 | 30.9 |

Table 6

$$\sum x = 217;$$
$$\sum y = 103.3;$$

$$\sum x^2 = 10{,}953;$$
$$\sum y^2 = 2{,}156.15;$$
$$\sum xy = 4{,}720.9;$$

$r_4 = r_i = 0.95075$ for any $i \in \{1, 2, \ldots, 10\}$.

Similarly if we replace in Pearson's formula (9) and also getting the same result equals to $r_p$:

3.5) Each $x_i$ by its deviation from x's median, and each $y_i$ by its deviation from y's median.

3.6) Each $x_i$ by its deviation from x's standard deviation, and each $y_i$ by its deviation from y's standard deviation.

3.7) Each $x_i$ by $x_i \pm a$ (where a is any number), and each $y_i$ by $y_i \pm b$ (where b is any number).

3.8) Each $x_i$ by $x_i * a$ (where a is any non-zero number and "*" is either division or multiplication), and each $y_i$ by $y_i * b$ (similarly for b and "*").

Since the cases 3.5 – 3.7 are similar to 3.1 - 3.4, let's consider two examples for the case 3.8:

3.8.1) Suppose each $x_i$ in the original example, Table 1, is divided by 5, while each $y_i$ is divided by 2.

Then:
$$\sum x = 71.4;$$
$$\sum y = 62.15;$$
$$\sum x^2 = 759.56;$$
$$\sum y^2 = 658.528;$$
$$\sum xy = 691.68;$$
$$r_{(\text{division, division})} = 0.95075.$$

3.8.2) Now, let's still divide each $x_i$ in Table 1 by 5, but this time multiply each $y_i$ with 2.

Then:
$$\sum x = 71.4;$$
$$\sum y = 248.6;$$

$$\sum x^2 = 759.56;$$
$$\sum y^2 = 10{,}536.4;$$
$$\sum xy = 2{,}766.72;$$
$$r_{(\text{division, multiplication})} = 0.95075.$$

So, again these results coincide with Pearson's.

More interesting alternative correlation coefficients [and given different results from Pearson's and Spearman's] are obtained by doing:

**A mixture of Pearson's and Spearman's correlation coefficients.**

4.1 We only replace $x_i$ by its rank among x's, while $y_i$ remains unchanged:

| x rank | 1 | 2 | 3 | 4 | 5 | 6 | 7 | 8 | 9 | 10 |
|---|---|---|---|---|---|---|---|---|---|---|
| y | 2.5 | 1.1 | 6.3 | 2.1 | 2.9 | 15.3 | 20.7 | 18.4 | 22 | 33 |

Table 7

$$\sum x = 55;$$
$$\sum y = 124.3;$$
$$\sum x^2 = 385;$$
$$\sum y^2 = 2{,}634.11;$$
$$\sum xy = 958.4;$$
$$r_{s,p} = 0.91661 \in [0.90303,\ 0.95075].$$

4.2. Similarly, as above, let's only replace $y_i$ by its rank among y's, while $x_i$ remains unchanged.

| x | 6 | 7 | 12 | 14 | 23 | 41 | 53 | 60 | 69 | 72 |
|---|---|---|---|---|---|---|---|---|---|---|
| y rank | 3 | 1 | 5 | 2 | 4 | 6 | 8 | 7 | 9 | 10 |

Table 8

$$\sum x = 357;$$
$$\sum y = 55;$$
$$\sum x^2 = 18{,}989;$$
$$\sum y^2 = 385;$$
$$\sum xy = 2{,}636;$$
$$r_{p,s} = 0.93698 \in [0.90303,\ 0.95075].$$

Both mixture correlation coefficients give different results from Pearson's and Spearman's, actually they are in between.

Conclusion:

In the samples where the rank in a discrete variable counts more than the variable values, this mixture of correlation coefficients brings better results than Pearson's or Spearman's.

Reference:

Jay Devore, Roxy Peck, "Introductory Statistics", second edition, West Publ. Co., 1994.